\documentclass{article}
\usepackage{amssymb,amsmath,latexsym}

\newcommand{\xdim}[3]{{}_{#2}^{#3}\,\!#1}
\newcommand{\xdimz}[2]{{}_{#1}^{#2}\,\!\mathbf{0}}
\newcommand{\sdim}[2]{{}^{#2}#1}
\newcommand{\sdimz}[1]{{}^{#1}\,\!\mathbf{0}}
\newcommand{\id}[1]{\sdim{\mathbf{I}}{#1}}
\newcommand{\vsp}[1]{\rule{0pt}{#1}}
\newcommand{\GFtwo}{\textsc{gf}{\small(2)}}

\newtheorem{thm}{Theorem}[section]
\newtheorem{definition}[thm]{Definition}
\newtheorem{lemma}[thm]{Lemma}

\newcommand{\finproof}{\quad\raisebox{-.225ex}{$\blacksquare$}}
\begin{document}
\title{Generating Large Non-Singular Matrices over an Arbitrary Field
					with Blocks of Full Rank}
\author{James Xiao \and Yongxin Zhou\thanks{yongxin.zhou@cloakware.com}
}
\date{June 27, 2002}

\maketitle
\begin{abstract}
\noindent
This note describes a technique for generating
large non-singular matrices with blocks of full rank.
Our motivation to construct such matrices 
arises in the white-box implementation of cryptographic algorithms with S-boxes.

\end{abstract}


\section{Introduction and Notation}

This note describes a technique for generating
large non-singular matrices with blocks of full rank.
One motivation is the following.
For ciphers such as \textsc{aes}\cite{DAERIJ01}, \textsc{des}\cite{FIPS46}, and
\textsc{cast}\cite{ADTA93}
involving linear transformations and substitution boxes (S-boxes),
white-box cryptographic implementations\cite{CEJVO02}
attempt to hide linear transformations in the non-linear S-box
lookups by blocking the matrices for the linear transformations,
and then non-linearly encoding the matrix operations
by converting the blocks into substitution boxes (S-boxes)
with arbitrary bijective input and output encodings.
Security considerations dictate that the matrices be hard to discover
from the S-boxes. A bijective S-box leaks no information if its input and
output codings are unknown and arbitrary, whereas a lossy S-box leaks
information: distinct encoded inputs map to the same encoded
output, reducing the search space for encodings. This in turn means that blocks
of reduced rank should be avoided.

We now introduce our notation.
Let $\xdim{M}{m}{n}$ denote an \mbox{$n \times m$}
matrix $M$ over field $F$;%
\footnote{We refer to this as $m$ inputs (columns) and $n$ outputs (rows),
because multiplying $\xdim{M}{m}{n}\,\xdim{X}{1}{m}$ yields a
vector $\xdim{Y}{1}{n}$.}
$\sdim{M}{n}$ is short for $\xdim{M}{n}{n}$.
$\id{n}$ denotes an \mbox{$n \times n$} identity matrix.
$\xdimz{m}{n}$ denotes an \mbox{$n \times m$} zero matrix;
$\sdimz{n}$ is short for $\xdimz{n}{n}$.
As usual, $m_{i,j}$ denotes the matrix $M$ element in \textit{row i} and
\textit{column j}.

A matrix may be \emph{blocked} into submatrices. For example,
\[
 \begin{pmatrix}
        \,              \xdim{A}{a}{b}  &       \xdim{B}{c}{b}  \\
        \vsp{12pt}\,    \xdim{C}{a}{d}  &       \xdim{D}{c}{d}
 \end{pmatrix}
\]
is a blocked matrix with blocks $A, B, C, D$, each of which is itself a
matrix. Horizontally adjacent
blocks must have the same number of rows, and vertically adjacent blocks must
have the same number of columns.

Where a matrix $\xdim{M}{m}{n}$ is blocked, and all of the blocks are square and
have the same dimensions
$p \times p$ with $p \vert m$ and $p \vert n$,
we use the notation $\xdim{M}{m}{n} [\sdim{B}{p}]$ to denote an $n\times m$
matrix $M$ with
$\frac{mn}{p^2}$ blocks.  Here $B_{i,j}$ denotes the block in row $i$
and column $j$ of blocks.

For convenience we give the following definition:

\begin{definition}
If all the blocks $B_{i,j}$ in a block matrix  $\xdim{M}{m}{n} [\sdim{B}{p}]$
are invertible, matrix
$M$ is called an $(m, n, p)$  \emph{block invertible matrix}. Furthermore, if
$m=n$, and $M$ is invertible then
$M$ is called an $(m,p)$ \emph{block invertible square matrix}.
\end{definition}

In this note we describe a way to create
a block invertible square matrix $\xdim{M}{n}{n} [\sdim{B}{p}]$
for $p$ and $n$ natural numbers where $p \vert n$ and $p > 1$.
One known technique involves the Kronecker product,
or tensor product of matrices\cite{Jacobson85}.
If we can find an invertible
matrix $\sdim{A}{p}$ such that all entries $a_{i,j}$ are not $0$ in field
$F$, its tensor product $ A \otimes B $
with another invertible matrix $\sdim{B}{p}$
is a $(p^2,p)$ block invertible square matrix.
However, this approach fails
for cases where the matrix $A$ does not exist --- for example,
when constructing $(2^t, 2)$ block invertible matrices over \GFtwo.
We provide a method of constructing block invertible matrices over any field.

\section{Preliminary Result}

First we prove the following result.

\begin{lemma}\label{lemma:id}
Let $p$ and $r$ be two integers with $ p > 1$ and $p \geq r \geq 0 $. Then there
exists a matrix $\sdim{A}{p}$ such that
\begin{gather*}
T ~=~
\begin{pmatrix} \id{r} &  & \\  &  \sdimz{p-r} \end{pmatrix}
 + A
\end{gather*}
is an invertible matrix over field $F$.
\end{lemma}

\noindent
\emph{Proof:}
We construct a matrix $\sdim{A}{p}$ such that $T$ is invertible, where
\begin{gather*}
T =
\begin{pmatrix}
\id{r} & \sdimz{p-r} \\ \sdimz{p-r} & \sdimz{p-r} \end{pmatrix}
+A ~.
\end{gather*}

\noindent
Case 1:
If $r$ is even, define
\begin{gather*}
A = \begin{pmatrix} 0 & 1  &    &    &    &    &    &     \\
                    1 & 1  &    &    &    &    &    &     \\
                      &    & 0  &  1 &    &    &    &     \\
                      &    & 1  &  1 &    &    &    &     \\
                      &    &    &   \dots & \dots   &    &     \\
                      &    &    &    &   0 & 1  &  &      \\
                      &    &    &    &   1 & 1  &  &       \\
                      &    &    &    &     &    &  \id{p-r}
\end{pmatrix}
\end{gather*}
with the $\tfrac{r}{2}$ invertible matrix
$\begin{pmatrix}0 & 1\\ 1 & 1 \end{pmatrix}$
on the diagonal. Therefore $A$ is invertible.
Since
\begin{gather*}
\id{2}+
\begin{pmatrix}0 & 1\\ 1 & 1 \end{pmatrix} =
\begin{pmatrix}1 & 1\\ 1 & 2 \end{pmatrix}
\end{gather*}
has a determinant of $1$, $T$ is invertible as each diagonal
block is invertible:

\begin{gather*}
T = \begin{pmatrix}\id{2}  &   &    &    &    &    &    &     \\
                      & \id{2}  &    &    &    &    &    &     \\
                      &    &   \dots & \dots   &    &     \\
                      &    &    &   \id{2} &   &  &       \\
                      &    &    &    &    &   \sdimz{p-r}
\end{pmatrix}
+ \begin{pmatrix} 0 & 1  &    &    &    &    &    &     \\
                    1 & 1  &    &    &    &    &    &     \\
                      &    & 0  &  1 &    &    &    &     \\
                      &    & 1  &  1 &    &    &    &     \\
                      &    &    &   \dots & \dots   &    &     \\
                      &    &    &    &   0 & 1  &  &      \\
                      &    &    &    &   1 & 1  &  &       \\
                      &    &    &    &     &    &  \id{p-r}
\end{pmatrix}
\end{gather*}
\begin{gather*}
=\begin{pmatrix}    1 & 1  &    &    &    &    &    &     \\
                    1 & 2  &    &    &    &    &    &     \\
                      &    & 1  &  1 &    &    &    &     \\
                      &    & 1  &  2 &    &    &    &     \\
                      &    &    &   \dots & \dots   &    &     \\
                      &    &    &    &   1 & 1  &  &      \\
                      &    &    &    &   1 & 2  &  &       \\
                      &    &    &    &     &    &  \id{p-r}
\end{pmatrix}.
\end{gather*}

\noindent
Case 2:
$r$ is odd.
For $r =1$, $A$ can be defined as
\begin{gather*}
A = \begin{pmatrix} 1 & 1  &     &     \\
                    1 & 0  &     &     \\
                      &    &  \id{p-2}
\end{pmatrix}
\end{gather*}

\noindent
Since  $p>1$, we have $p-2\geq 0$. Note that $A$ is invertible, and
\begin{gather*}
T = \begin{pmatrix}1 & 0 &   &           \\
                   0 & 0 &   &           \\
                     &   &  \sdimz{p-2} \end{pmatrix}
+ A =
\begin{pmatrix} 0 & 1  &   &         \\
                1 & 0  &   &         \\
                  &    &  \id{p-2} \end{pmatrix}
\end{gather*}
is also invertible.

For $r$ odd and $r>1$, let $r = 2n + 3$ where $n\geq 0$.
Now define $A$ to be

\begin{gather*}
A = \begin{pmatrix} 1 & 1  & 1  &    &    &    &    &     \\
                    1 & 1  & 0  &    &    &    &    &     \\
                    1 & 0  & 0  &    &    &    &    &     \\
                      &    &    &  0 &  1 &    &    &     \\
                      &    &    &  1 &  1 &    &    &     \\
                      &    &    &   \dots & \dots    &    &     \\
                      &    &    &    &   0 & 1 &  &      \\
                      &    &    &    &   1 & 1 &  &       \\
                      &    &    &    &     &   &  \id{p-2n}
\end{pmatrix}
\end{gather*}

\noindent
Note that
\begin{gather*}
\id{3} + \begin{pmatrix} 1 & 1  & 1  \\
                1 & 1  & 0  \\
                1 & 0  & 0 \end{pmatrix}
=\begin{pmatrix} 0 & 1  & 1  \\
                1 & 0  & 0  \\
                1 & 0  & 1 \end{pmatrix}
\end{gather*}
which is invertible; and using the same argument as above,
both $A$ and $T$ are invertible.\finproof

\section{Constructing a block invertible square matrix}

Before proceeding, we recall an elementary result used in the proof of our main result.

\begin{lemma}\emph{From Paley and Weichsel\cite{Paley72}:}\label{lemma:id1}
For a given square matrix $\sdim{M}{n}$ of rank
$r \leq n$ over field $F$, there exist invertible matrices $\sdim{P}{n}$
and $\sdim{Q}{n}$ such that
\begin{gather*}
M = P
\begin{pmatrix} \id{r} &  &\\   & \sdimz{n-r} & \end{pmatrix}
 Q.
\end{gather*}
\end{lemma}

\begin{thm}[Main result]
For any field $F$, and for any positive integers $n$ and $p$ such that
$n \geq p$ and $p\vert n$, there exists an $(n, p)$ block invertible square
matrix.
\end{thm}

\noindent
\emph{Proof:}
We construct the matrix inductively.
For the first step we find an invertible square matrix $\sdim{M}{p}$ over $F$.
Note that there are infinitely many $p \times p$ invertible matrices over
infinite field $F$ and there are 
$$\prod_{i=0}^{p-1}(q^{p}-q^{i})$$ 
invertible matrices over finite field $F$ of order $q$\cite{Dickson58}.
These facts
grant us a variety of choices of $\sdim{M}{p}$.
This $\sdim{M}{p}$ is a $(p, p)$ block invertible square matrix.

Now suppose we have found a $(t, p)$ block invertible square matrix $M$ with
$t \geq  p$ and $p \vert t$.
The third step is
to construct a $(t+p, p)$ block invertible square matrix.

It is not hard to see that there exists a $(p, t, p)$ block invertible matrix
$X$ and a $(t, p, p)$ block invertible matrix $Y$.
In fact $X$ and $Y$ can be constructed from $M$ because $M$ is a $(t, p)$
block invertible square matrix. 

Let $\sdim{W}{p}$ be a matrix over $F$. Observing the following matrix equation:

\begin{gather*}
\begin{pmatrix} M & 0 \\  X & W \end{pmatrix} \cdot
\begin{pmatrix} I & M^{-1}Y \\  0 & I \end{pmatrix} =
\begin{pmatrix} M & Y \\  X & XM^{-1}Y + W \end{pmatrix}
\end{gather*}

\noindent
we claim that if we can find a $p\times p$ invertible matrix $W$ such that
$XM^{-1}Y + W$ is invertible, then matrix
\begin{gather*}
N = \begin{pmatrix} M & Y \\  X & XM^{-1}Y + W \end{pmatrix}
\end{gather*}
is a $(t+p, p)$ invertible square matrix. In fact, if $W$ is invertible, the
left-side of the matrix equation implies $N$ is invertible. Following the
assumptions that
$M$, $X$, $Y$, and $ XM^{-1}Y+W$ are $(t,p)$, $(p,t, p)$, $(t, p, p)$
and $(p,p)$ block invertible matrices, respectively, by definition, $N$ is a
block invertible square matrix.
Such a matrix $W$ can be constructed in the following way.
  
>From Lemma \ref{lemma:id1},
for the $p\times p$ square matrix $XM^{-1}Y$,
there exist two invertible
matrices $\sdim{P}{p}$  and $\sdim{Q}{p}$ such that

\begin{gather*}
P(XM^{-1}Y)Q =
\begin{pmatrix}
\id{r} & \sdimz{p-r} \\ \sdimz{p-r} & \sdimz{p-r} \end{pmatrix}
\end{gather*}
where $r$ is the rank of $XM^{-1}Y$.
By Lemma \ref{lemma:id}, an invertible matrix $\sdim{A}{p}$ exists such that
\begin{gather*}
\begin{pmatrix}
\id{r} & \sdimz{p-r} \\ \sdimz{p-r} & \sdimz{p-r} \end{pmatrix}
+A
\end{gather*}
is invertible. Now we can define $W$ as $W = P^{-1}AQ^{-1}$.
Then the matrix
\begin{gather*}
 XM^{-1}Y + W = P^{-1}( P(XM^{-1}Y)Q + A)Q^{-1}
 =  P^{-1}(
\begin{pmatrix}
\id{r} & \sdimz{p-r} \\ \sdimz{p-r} & \sdimz{p-r} \end{pmatrix}
+A)Q^{-1}
\end{gather*}
is invertible, completing our construction.\finproof

\section{Example}

In the following example, for $F=$ \GFtwo,
we construct an $(n, 2)$ block invertible
square matrix for any even number $n>0$.
All blocks mentioned below are of dimension $2\times 2$.

Since any invertible $2\times 2$ matrix is a $(2,2)$ block invertible square
matrix, we
can safely assume that we already have a $(t, 2)$ block invertible square matrix
$\sdim{M}{t}$ for $t \geq 2$.
We construct a $(t+2, 2)$ block invertible square matrix
$\sdim{M'}{t+2}$ as follows.

\begin{enumerate}
\item Use a row of blocks of matrix $M$ to create a matrix $\xdim{X}{t}{2}$.
\item Use a column of blocks of matrix $M$ to create a matrix $\xdim{Y}{2}{t}$.
\item Get invertible matrices $\sdim{P}{2}$  and $\sdim{Q}{2}$ such that
      $P(XM^{-1}Y)Q = \begin{pmatrix} \id{r} &  \\ & \sdimz{2-r} \end{pmatrix}$.

\item Define matrix $\sdim{A}{2}$ as follows:
	\begin{enumerate}
           \item if $r = 0$, $A=\id{2}$;
           \item if $r = 1$, $A = \begin{pmatrix}1&1\\1&0\end{pmatrix}$;
           \item if $r = 2$, $A = \begin{pmatrix}0&1\\1&1\end{pmatrix}$.
	\end{enumerate}

\item Then $\begin{pmatrix}M & Y \\ X & XM^{-1}Y + P^{-1}AQ^{-1}\end{pmatrix}$
      is a $(t+2, 2)$ block invertible matrix.
\end{enumerate}

\noindent
Repeat for $(n-2)/2$ steps to obtain an $(n,2)$ block invertible square matrix.


\end{document}